\documentclass[12pt]{article}

\usepackage{geometry}
\usepackage{graphicx}
\usepackage{lipsum}
\usepackage[utf8]{inputenc}
\usepackage{amsmath,amscd}    
\usepackage{graphicx}   
\usepackage{verbatim}   
\usepackage{color}      
\usepackage{subfigure}  
\usepackage{hyperref}   
\usepackage{mathtools}
\usepackage{amssymb}
\definecolor{titlepagecolor}{cmyk}{1,.60,0,.40}
\definecolor{namecolor}{cmyk}{1,.50,0,.10}

\newcommand{\lbiprod}{{>\!\!\!\triangleleft\kern-.33em\cdot}}
\newcommand{\rbiprod}{{\cdot\kern-.33em\triangleright\!\!\!<}}

\begin{document}

\author{Niels Gresnigt}

\title{Knotted boundaries and braid only form of braided belts}
\maketitle
\begin{abstract}

The Helon model identifies Standard Model quarks and leptons with certain framed braids joined together at both ends by a connecting node (disk). These surfaces with boundary are called braided 3-belts (or simply belts). Twisting and braiding of ribbons composing braided 3-belts are interchangeable, and it was shown in the literature that any braided 3-belt can be written in a pure twist form, specified by a vector of three multiples of half integers $[a,b,c]$, a topological invariant.    

This paper identifies the set of braided 3-belts that can be written in a braid only form in which all twisting is eliminated. For these braids an algorithm to calculate the braid word is determined which allows the braid only word of every braided 3-belt to be written in a canonical form. It is furthermore demonstrated that the set of braided 3-belts do not form a group, due to a lack of isogeny.

The conditions under which the boundary of a braided 3-belt is a knot are determined, and a formula for the Jones polynomial for knotted boundaries is derived. Considering knotted boundaries makes it possible to relate the Helon model to a model of quarks and leptons in terms of quantum trefoil knots, understood as representation of the quantum group $SU_q(2)$. Associating representations of a quantum group to the boundary of braided belts provides a possible means of developing the gauge symmetries of interacting braided belts in future work.
\end{abstract}

\newpage
\tableofcontents

\section{Introduction}

In 2005, Bilson-Thompson proposed the Helon model in which one generation of leptons and quarks are identified as braids of three ribbons connected at the top and bottom via a disk \cite{Bilson-Thompson2005}. These framed braids are called braided 3-belts, and they are surfaces with boundary \cite{Bilson-Thompson2009,shepperd1962braids}. With the additional structure that each ribbon can be twisted clockwise or anticlockwise by $2\pi$ (interpreted physically as electric charge), and satisfying certain conditions, the simplest non-trivial braided belts map precisely to the first generation fermions of the Standard Model (SM). The original model has since been expanded into a complete scheme for the identification of the SM fermions and weak vector bosons for an unlimited series of generations \cite{Bilson-Thompson2009,Bilson-Thompson2012,Bilson-Thompson2008}. 

These topological representations of leptons and quarks may be embedded within a larger network of braided ribbons. Such a ribbon network is a generalization of a spin network, fundamental in Loop Quantum Gravity (LQG) \cite{Bilson-Thompson2005,Bilson-Thompson2007,Bilson-Thompson2009}. The embedding of framed braids into ribbon networks makes it possible to develop a unified theory of matter and spacetime in which both are emergent from the ribbon networks.

It was shown in \cite{Bilson-Thompson2009} that it is possible to deform a braided 3-belt in a manner that interchanges between the braiding of ribbons and the twisting of ribbons. In the same paper it was shown that it is always possible to deform the belt into a structure that carries only twisting (the braiding is trivial). This form is called the \textit{pure twist form} and has associated with it a \textit{pure twist word} $[a,b,c]$ with $a,b,c$ being multiples of half integers denoting the number of $\pm \pi$ twists carried by each ribbon. This pure twist word is a topological invariant and so every braided belt can be classified into an equivalence class specified by the pure twist word.

One might wonder if the converse is likewise possible. That is, can all the twists on the ribbons be systematically removed through deforming the belt such that one is left with a belt that is only braided? The first part of this paper addresses this question. Although not possible in general, this is possible for the framed 3-braids of the Helon model. More generally it is possible for any orientable braided 3-belt

There are situations where writing a braided 3-belt in terms of braiding only rather than twisting only is advantageous. Together with the observation that Clifford algebras and normed division algebras contain representations of the circular Artin braid groups \cite{kauffman2016braiding,gresnigt2017quantum2}, it was recently shown that by removing the twists on ribbons in the Helon model braids, they can be written in a form such that they coincide with the basis states of the minimal left ideals of the complex octonions \cite{gresnigt2018braids}, which have been shown to transform as a single generation of leptons and quarks under the unbroken symmetries $SU(3)_c$ and $U(1)_{em}$ \cite{furey2016standard,stoica2017standard}. It is precisely because all the twisting can be exchanged for braiding that this identification between Helon braids and basis states of the minimal ideals of the complex octonions is possible. This curious connection between two complementary models is worth investigating further.

Having established when it is possible to write a braided 3-belt in a braid only form, an algorithm for determining the braid word is then presented. Unlike the pure twist word $[a,b,c]$, the braid word is not a topological invariant. It is however always possible to rewrite the braid word in a canonical form. 

The boundary of a braided belt is a link. The reduced link, defined as the boundary link with any unknotted and unlinked loops removed is an invariant under the evolution algebra generated by the Pachner moves \cite{Bilson-Thompson2007}. These boundary links of braided belts are studied in the second part of this paper, and it is determined under what conditions the boundary link is a knot. A formula is then derived for the Jones polynomial of the knotted boundary in terms of the pure twist word $[a,b,c]$. The knotted boundaries of the simplest braided belts are then tabulated. Although a rigorous study of the spectrum of knots that can form the boundaries of braided belts is left for a later date, this initial work hints that the knots are always prime knots. 

By classifying braided 3-belts in terms of their boundary knots it becomes possible to find alternate mappings from braided belts to leptons and quarks. For example, a model by Finkelstein represents leptons and quarks as quantum trefoil knots based on representations of the quantum group $SU_q(2)$ \cite{Finkelstein2007,Finkelstein2010,finkelstein2013preon}. Because the trefoil knot can form the boundary of a braided belt, this knot model provides an alternate mapping from braided belts to SM fermions. Concurrently, considering these trefoil knots as the boundaries of braided 3-belts provides an embedding of Finkeltein's model into ($q$-deformed) spin-networks, and hence LQG. 

The outline of the paper is as follows. A brief overview of the Helon model is given in section 2 before the Artin braid groups and braided belts are reviewed in section 3. Section 4 shows that braid only form of a braided 3-belt exists if and only if the belt is orientable, and subsequently find an algorithm for determining the braid only word. Section 5 looks at when the boundary of a braided belt forms a knot and subsequently a formula for the Jones polynomial is derived. The knotted boundaries for the simlest cases are then tabulated. Finally, in section 6 a comparison between braided belts and the Helon model with an alternative model of leptons and quarks as trefoil knots is made. 

\section{The Helon model}\label{helonmodelsection}

\begin{figure}[h!]
\centering
\includegraphics[scale=0.30]{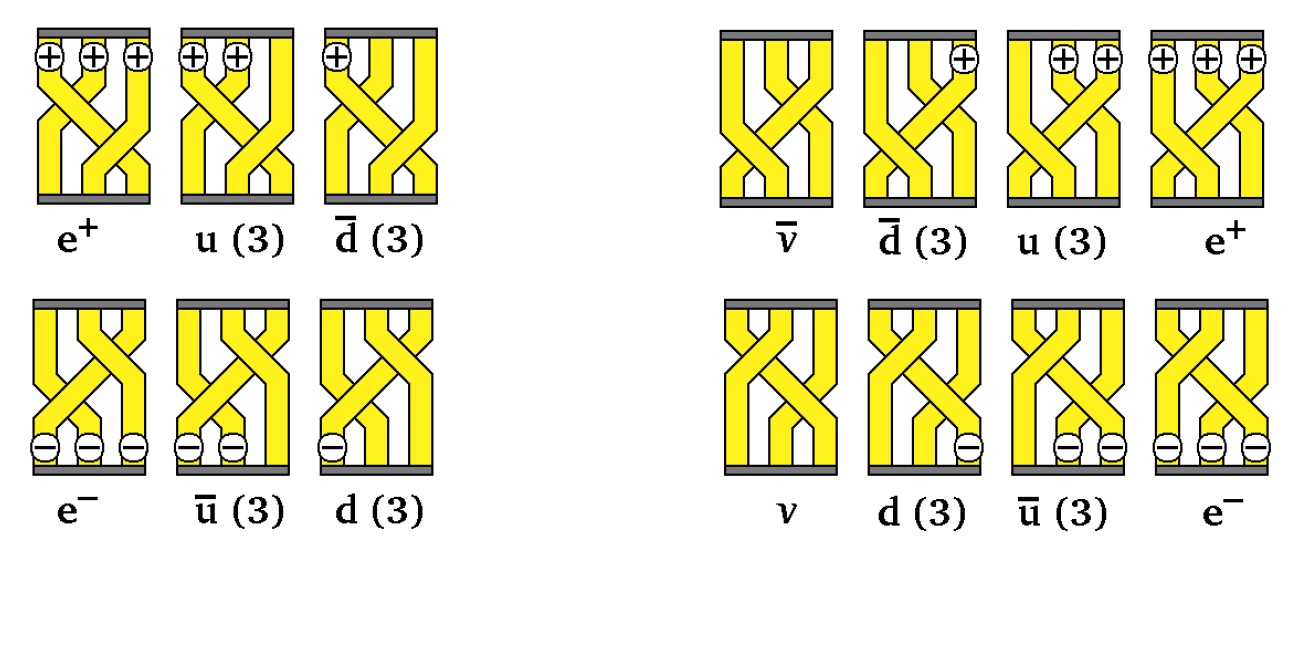}
\caption{The Helon model of Bilson-Thompson in which the first generation SM fermions are represented as braids of three (possibly twisted) ribbons. Used with permission. Source, \cite{Bilson-Thompson2005}.}
\label{helonmodel}
\end{figure}

The Helon model of Bilson-Thompson maps the simplest non-trivial braids consisting of three (twisted) ribbons and two crossings to the first generation of SM fermions. Quantized electric charges of particles are represented by integral twists of the ribbons of the braids, with a twist of $\pm 2\pi$ representing an electric charge of $\pm e/3$. The twist carrying ribbons, called helons, are combined into triplets by connecting the tops of three ribbons to each other and likewise for the bottoms of the ribbons. The color charges of quarks and gluons are accounted for by the permutations of twists on certain braids, and simple topological processes are identified with the electroweak interaction, the color interaction, and conservation laws. The lack of twist on the neutrino braids means they only come in one handedness. The representation of first generation SM fermions in terms of braids is shown in Figure \ref{helonmodel}.

Braided belts may be embedded within a larger network of braided ribbons. Such a ribbon network is a generalization of a spin network. Spin networks are, in the simplest case, trivalent graphs with an irreducible representation of $SU(2)$ and an intertwiner associated to every link and vertex respectively. Instead labelling the links by representations of the quantum group $SU_q(2)$ introduces a nonzero cosmological constant into the theory and has the effect that the links become framed \cite{turaev1992state,Bilson-Thompson2007}.

Within ribbon networks, braided belts correspond to local noiseless subsystems which have been shown to exist in background independent theories where the microscopic quantum states are defined in terms of the embedding of a framed, or ribbon, graph in a three manifold, and in which the allowed evolution moves are the standard local exchange and expansion moves (Pachner moves). Such noiseless subsystems are given by braided sets of $n$ edges joined at both ends by a set of connected nodes, or equivalently, a disk. The embedding into a ribbon network is possible by connecting (at least) one of the disks to the rest of the ribbon graph. What the Helon model shows is that the simplest emergent local structures of such theories, when $n=3$, match precisely the first generation leptons and quarks. It has furthermore been shown that an embedding into tetravalent networks is likewise possible \cite{Bilson-Thompson2012}. 

Discrete symmetries have already been studied in the Helon model and may be defined on the braid in such a way that performing all three in any order leaves the braid unchanged \cite{he2008c}. Dynamics and interactions of braids have been studied in terms of evolution moves on trivalent and tetravalent spin network. In trivalent spin networks, the braid excitations are too strongly conserved meaning a lack of creation and annihilation. However Smolin and Wan have shown that in the tetravalent case, braid excitations can propagate and interact under the dual Pachner moves \cite{smolin2008propagation} . The dynamics of these braid excitations can be described by an effective theory based on Feynman diagrams \cite{wan2009effective}. 

\section{The Artin braid groups and braided belts}

\subsection{The Artin braid groups}

The Artin braid group on $n$ strands is denoted by $B_n$ and is generated by elementary braids $\left\lbrace \sigma_1,...,\sigma_{n-1}\right\rbrace$ subject to the relations
\begin{eqnarray}\label{braidrelations}
\sigma_i \sigma_j=\sigma_j\sigma_i,\;\textrm{whenever}\; \vert i-j \vert > 1,\\
\sigma_i\sigma_{i+1}\sigma_i=\sigma_{i+1}\sigma_i \sigma_{i+1},\;\textrm{for}\; i=1,....,n-2.\label{braidrelations2}
\end{eqnarray}
The braid groups $B_n$ are an extension of the symmetric groups $S_n$ with the condition that the square of each generator being equal to one lifted. The inverse of a braid is the unique braid which unbraids the first braid when the two are composed. Drawing braids vertically as in figure \ref{helonmodel}, the braid inverse corresponds to a vertical reflection of the braid. The braid inverse is an anti-automorphism so that, for example, $(\sigma_3\sigma_1\sigma_2^{-1})^{-1}=\sigma_2\sigma_1^{-1}\sigma_3^{-1}$. Two braids are said to equivalent (or isotopic) when they can be continuously deformed into each other inside $\mathbb{R}^3                                 $.

Replacing the strands by ribbons provides a framing of the braid. In the framed braid group, each strand is thickened to a ribbon with the additional structure that now each ribbon may be twisted. Thus in addition to the braid generators $\sigma_1,...,\sigma_{n-1}$ of $B_n$, the framed braid group has the additional twist operators $t_1,...,t_n$. The framed braid group of $B_n$ is then defined by relations (\ref{braidrelations}, \ref{braidrelations2}) and additional relations
\begin{eqnarray}
t_it_j &=&t_jt_i,\;\textrm{for}\;\textrm{all}\;i,j,\\
\sigma_it_j&=&t_{\sigma_i(j)}\sigma_i,
\end{eqnarray}
where $\sigma_i(j)$ denotes the permutation induced on $j$ by $\sigma_i$. For example $\sigma_1(2)=1$ and $\sigma_1(3)=3$.

A circular braid on $n$ strands has $n$ strands attached to the outer edges of two circles which lie in parallel planes in $\mathbb{R}^3$. The circular braid group on $n$ strands is denoted $B_n^c$, and has $n$ generators. Specifically, the circular braid group $B_3^c$ has three generators $\sigma_1,\sigma_2,\sigma_3$, satisfying
\begin{eqnarray}
\sigma_1\sigma_2\sigma_1=\sigma_2\sigma_1\sigma_2,\;\;\sigma_2\sigma_3\sigma_2=\sigma_3\sigma_2\sigma_3,\;\;\sigma_3\sigma_1\sigma_3=\sigma_1\sigma_3\sigma_1.
\end{eqnarray}
In the Helon model, the ribbons are connected at both ends to framed trivalent nodes, which are disks. This is equivalent to connecting the ribbons to two parallel disks, forming a surface with boundary. The appropriate braid group is therefore not just $B_3$ but rather the circular braid group $B_3^c$. Unlike in $B_3$, in $B_3^c$ each strand or ribbon is treated symmetrically. In much of what follows we will be using this fact to write expressions as symmetrically as possible in terms of not only $\sigma_1$ and $\sigma_2$, but rather $\sigma_1,\sigma_2,\sigma_3$. 

We can multiply two framed braids together by first joining the bottom of the ribbons of the first braid to the tops of the ribbon of the second braid and then sliding (isotoping) the twists from each component braid downward. Doing so, the twists carried by the first braid will get permuted by the second braid. We will write the resulting product in standard form with the braiding first followed by the twisting. We can then write the composition law as 
\begin{eqnarray}\label{braidproduct}
\nonumber (\Lambda_1,[a_1,a_2,a_3])(\Lambda_2,[b_1,b_2,b_3])&=&(\Lambda_1\Lambda_2, P_{\Lambda_2}([a_1,a_2,a_3])[b_1,b_2,b_3]),\\
\nonumber&=&(\Lambda_1\Lambda_2, [a_{\pi(\Lambda_2)(1)},a_{\pi(\Lambda_2)(2)},a_{\pi(\Lambda_2)(3)}])[b_1,b_2,b_3],\\
&=&(\Lambda_1\Lambda_2, ([a_{\pi(\Lambda_2)(1)}+b_1,a_{\pi(\Lambda_2)(2)}+b_2,a_{\pi(\Lambda_2)(3)}+b_3])
\end{eqnarray}
where $\Lambda_1$ and $\Lambda_2$ are two braid words, $P_{\Lambda_i}$ is the permutation induced on $[a,b,c]$ by the braid word $\Lambda$, $[a_1,a_2,a_3][b_1,b_2,b_3]=[a_1+b_1, a_2+b_2, a_3+b_3]$, and $\pi:B_3^c\rightarrow S_3$ with $\pi(\sigma_1)=(12)$, $\pi(\sigma_2)=(23)$, $\pi(\sigma_3)=(31)$. 

\subsection{Braided belts}
\begin{figure}[h!]
\centering
   \includegraphics[scale=0.3]{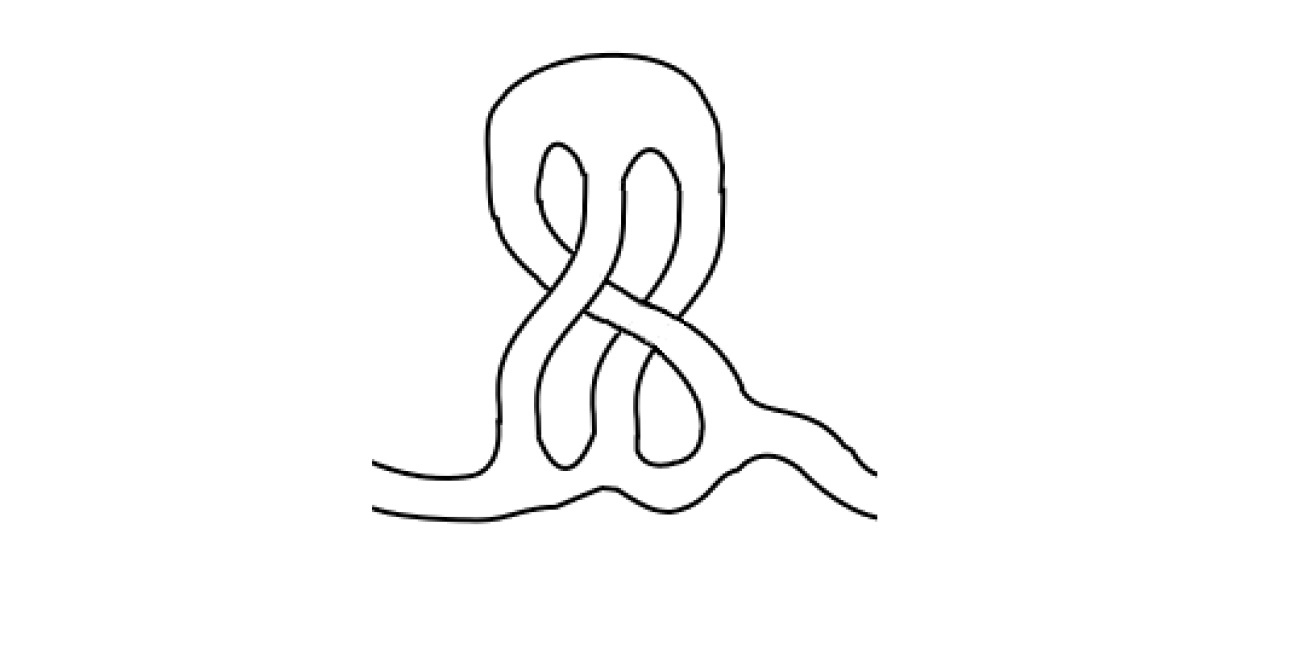}
\caption{The embedding of a braided belt into a ribbon network. Used with permission. Source, \cite{Bilson-Thompson2007}.}
\label{embedding}
\end{figure}

Connecting three braided ribbons together by a disk at both ends makes a braided belt. A braided belt is a surface with boundary composed of the union of the disks and ribbons. A braided 3-belt may be embedded within a larger trivalent ribbon network by connecting one end of the braid to the network whilst the other end of the braid remains connected to a disk. The embedding of a braided 3-belt into a ribbon network is shown in Figure \ref{embedding}. 

\begin{figure}[h!]
\centering
   \includegraphics[scale=0.15]{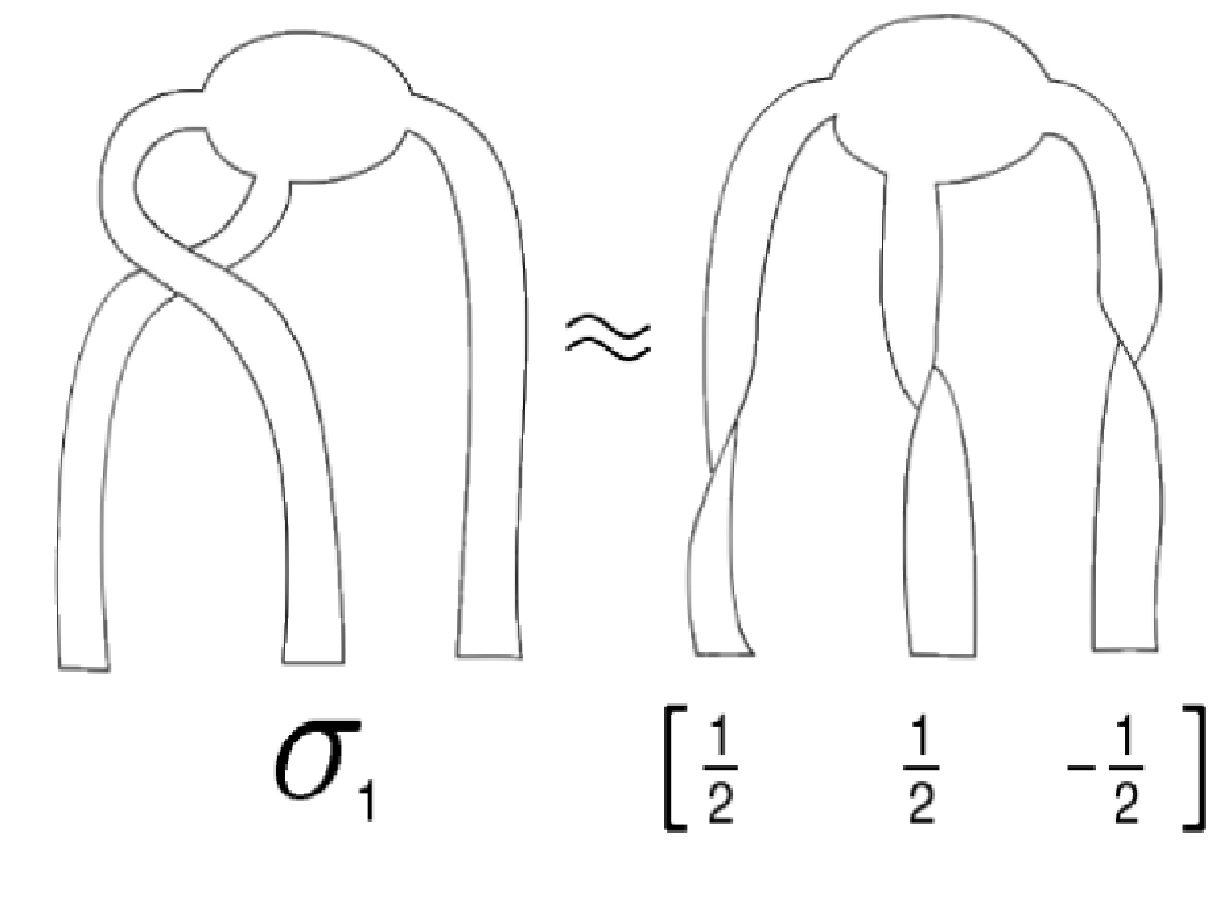}
\caption{The twisting and braiding in framed 3-braids connected at both ends to a disk is interchangeable. Source, \cite{Bilson-Thompson2007}.}
\label{braidtotwist2}
\end{figure}

A crucial point in what follows is that the twisting and the braiding of braided 3-belts are not individually topologically invariant. Through suitably manipulating one of the disks to which a braid is connected it is possible to interchange between braiding and twisting \cite{Bilson-Thompson2008}. This is achieved by flipping the disk at one end over or passing it through the ribbons. In this way the braiding may be reduced and exchanged for twisting. In Figure \ref{braidtotwist2}, it is shown how the braiding induced by the generator $\sigma_1$ may be exchanged for twisting\footnote{Not shown explicitly in the figure are the disks are the bottom of the two braids. These must be included for the interchanging of braiding with twisting to work. Without these disks it would be possible to simply slide the first two ribbons in the first braid past each other without introducing any twisting at all.}. 

\subsection{Pure twist form of a braided belt}

It was shown in \cite{Bilson-Thompson2009} that for braided 3-belts it is always possible to iteratively remove all the braiding of the ribbons. The result is an unbraided structure of three ribbons that each carry half integer number of half twists. This is only possible for the case of a 3-belt and does not hold in general for a braided belt composed of a higher number of ribbons. One may then write the braided 3-belt in \textit{pure twist form} as a vector $[a,b,c]$ where $a,b,c$ are three integer multiples of $1/2$. We follow \cite{Bilson-Thompson2009} and refer to the vector $[a,b,c]$ as the \textit{pure twist word}.

The braid generators of the circular Artin braid group $B^c_3$ can be written as twist vectors as follows:
\begin{eqnarray}\label{braidgenerators}
\nonumber\sigma_1\rightarrow \left[ \frac{1}{2},\frac{1}{2},-\frac{1}{2}\right] ,\qquad \sigma^{-1}_1\rightarrow \left[ -\frac{1}{2},-\frac{1}{2},\frac{1}{2}\right] ,\\
\nonumber\sigma_2\rightarrow \left[ -\frac{1}{2},\frac{1}{2},\frac{1}{2}\right] ,\qquad\sigma^{-1}_2\rightarrow \left[ \frac{1}{2},-\frac{1}{2},-\frac{1}{2}\right],\\
\sigma_3\rightarrow \left[ \frac{1}{2},-\frac{1}{2},\frac{1}{2}\right] ,\qquad\sigma^{-1}_3\rightarrow \left[ -\frac{1}{2}, \frac{1}{2},-\frac{1}{2}\right] .
\end{eqnarray}

One advantage of writing a braid in this pure twist form is that the twist vector is a topological invariant. Braids with the same pure twist word are isotopic. This means that each 3-belt belongs to an equivalence class specified by the triplet $[a,b,c]$. A second feature of the pure braid word is that it does not mix integers and half integers if and only if the surface with boundary corresponding to a braided 3-belt is orientable. In such cases we say the the braided 3-belt is orientable.

Because twisting and braiding do not in general commute, one must take care when turning a general braid word into a twist word as each braiding permutes two of the ribbons. Thus the twists on ribbons get permuted according to the permutation $P_B$ associated with the braid $B$. We write the  twist vector permuted by the braid $B$ as $P_{B}([a,b,c])$. For example, we have $P_{\sigma_1}([a,b,c])=[b,a,c]$ and $P_{\sigma^{-1}_2}([a,b,c])=P_{\sigma_2}([a,b,c])=[a,c,b]$. 

\subsection{The product of braided belts and the lack of group structure}

As with framed braids, and described above in equation (\ref{braidproduct}), one composes braided belts by joining the bottom of the ribbons of the first belt to the tops of the ribbons of the second. Doing so requires that the intervening nodes (disks) be removed. However, although framed braids generate a group, we will show here that the set of braided belts do not. This is because the composition described above is not defined up to isogeny. That is, if $b_1$ and $b'_1$ are two isotopic braided belts and $b_2$ and $b'_2$ are likewise two isotopic braided belts, then is does not follow that the composition $b_1b_2$ is isotopic to $b'_1b'_2$. For example, note that
\begin{eqnarray}
\sigma_2\sigma_1[0,0,0]=\sigma_3\sigma_1[0,0,0]\rightarrow[0,0,1],
\end{eqnarray}
and so $\sigma_2\sigma_1$ and $\sigma_3\sigma_1$ are isotopic braided belts. Likewise,
\begin{eqnarray}
\sigma_3\sigma_2[0,0,0]=\sigma_1\sigma_2[0,0,0]\rightarrow[1,0,0],
\end{eqnarray}
so that $\sigma_3\sigma_2$ and $\sigma_1\sigma_2$ are isotopic braided belts. However,
\begin{eqnarray}
\sigma_1\sigma_2\sigma_2\sigma_1[0,0,0]\rightarrow [2,0,0],
\end{eqnarray}
whereas
\begin{eqnarray}
\sigma_3\sigma_2\sigma_2\sigma_1[0,0,0]\rightarrow [1,1,0].
\end{eqnarray}
Hence the isogeny of the composition law fails to hold, and braided belts do not form a group. The original of this breakdown of isogeny is the different permutations induced by $\sigma_1\sigma_2$ and $\sigma_3\sigma_2$. In general we have
\begin{eqnarray}
\sigma_1\sigma_2[a,b,c]\rightarrow[c+1,a,b],\qquad \sigma_3\sigma_2[a,b,c]\rightarrow[b+1, c,a],
\end{eqnarray}
which in general are not isotopic. A group structure is however recovered for the subset of braided 3-belts that do not permute their composing ribbons. (Framed) braids that to do not permute their (ribbons) strands make up the kernel of the homomorphism $\pi:B_n^c\rightarrow S_n$. More will be said about this in the next section, where the kernel is used to solve the word problem for braided 3-belts. We note that the lack of group structure for braided belts was not discussed in the earlier work on braided belts \cite{Bilson-Thompson2009}.

\section{Braid only form of a braided 3-belt}

Several methods for finding the pure twist form of a 3-belt are outlined in \cite{Bilson-Thompson2009}. The present paper complements that work and considers under what conditions a braided 3-belt may be written in a \textit{braid only form} where instead of the braiding, the twists on the ribbons are made trivial. An algorithm is then outlines for determining the braid word. Unlike the pure twist vector, the braid word is not a topological invariant, and is thus not unique. Algorithms however exist that can be used to write any braid in a normal form. One example is Artin's combing algorithm which solves the word problem and allows any braid to be written in a normal combed form \cite{artin1947theory}. 

As noted earlier, the pure twist word of an orientable braided 3-belt does not mix integers and half integers. Multiplying a braid $B$ on the left (or right) by any braid generator preserves the orientability (or non-orientability) of the $B$. This is readily seen by considering equations (\ref{braidgenerators}). The twist vectors for the braid generators do no mix integers and half integers and it follows that composing them with any pure twist word does not induce mixing or unmixing of integers and half integers (regardless of any permutation that is introduced). Because the twist word $[a,b,c]$ is a topological invariant, a pure twist word is trivial (isotopic to the trivial braided 3-belt) if and only if $a=b=c=0$. Because this trivial braided 3-belt $[0,0,0]$ is orientable, the action of any number of braid generators on this identity will preserve its orientablity. In other words, braided 3-belts that can be written in a braid only form in which the twisting is trivial are necessarily orientable. 

The converse likewise holds. That is, any orientable braided 3-belt has an associated braid only form. It is readily verified that
\begin{eqnarray}
\sigma_1\sigma_2[0,0,0]= [1,0,0],\qquad \sigma_2\sigma_3[0,0,0]=[0,1,0],\qquad \sigma_3\sigma_1[0,0,0]=[0,0,1].
\end{eqnarray}
The negatives can be obtained by replacing the braid generators by their inverses (keeping the order the same)\footnote{Furthermore, $(\sigma_1\sigma_2)^2=[1,1,0]$, $(\sigma_2\sigma_3)^2=[0,1,1]$,  $(\sigma_3\sigma_1)^2=[1,0,1]$, and $(\sigma_1\sigma_2)^3=(\sigma_2\sigma_3)^3=(\sigma_3\sigma_1)^3=[1,1,1]$.}. Notice that using all three generators of the circular braid group $B_3^c$ highlights the symmetry in this situation. It follows then that any pure twist vector $[a,b,c]$ with integer $a,b,c$ can be written in braid only form via a clever combination of the above operators. We say clever, because one needs to keep track of the permutations induced. For example, $(\sigma_2\sigma_3)[0,0,0]=[0,1,0]$ however $[0,2,0]\neq (\sigma_2\sigma_3)^2[0,0,0]$ but rather $(\sigma_2\sigma_3)(\sigma_1\sigma_2)[0,0,0]=[0,2,0]$. More generally,
\begin{eqnarray}\label{eq7}
\nonumber\sigma_1\sigma_2[a,b,c]=[c+1,a,b],\quad \sigma_2\sigma_3[a,b,c]=[c,a+1,b],\quad \sigma_3\sigma_1[a,b,c]=[c,a,b+1].
\end{eqnarray}
It remains only to show that when $a,b,c$ are half integers instead, one can likewise find a pure braid word. To do so we consider $[a-\frac{1}{2},b-\frac{1}{2},c-\frac{1}{2}]$ for which the existence of a braid only form has already been established. It is readily checked that $\sigma_1\sigma_2\sigma_1[a,b,c]=[c+\frac{1}{2},b+\frac{1}{2}, a+\frac{1}{2}]$. By writing the pure braid word for  $[c-\frac{1}{2},b-\frac{1}{2},a-\frac{1}{2}]$ and then acting on it with $\sigma_1\sigma_2\sigma_1$, from the left, a braid only word for $[a,b,c]$ can be obtained.

We conclude that a braided 3-belt can be written in braid only  form if and only if the 3-belt is orientable. From now on only orientable belts will be considered. Having shown that all orientable braided 3-belts can be written in pure braid from, we now present an algorithm for determining the pure braid word for an orientable 3-belt. 

\subsection{The word problem for pure 3-belts}

The braid word problem is to decide whether two different braid words represent the same braid. The word problem has been solved, first by Artin in 1947 \cite{artin1947theory}. Artin's algorithm makes use of the group theoretic properties of the kernel of the homomorphism $\pi:B_n\rightarrow S_n$ to put the braid into a normal form, called a combed braid. 

The braid only form of an orientable 3-belt is not a topological invariant since it is not unique. It is then desirable to write a braided belt in a normal form. This is equivalent to solving the word problem for braided 3-belts. This turns out to be surprisingly simple. The reader is reminded that what follows is true only for braided belts consisting of three ribbons, not for braids, of braided belts in general.

The kernel of the homomorphism $\pi : B_3^c\rightarrow S_3$ is the pure braid group $PB_3^c$ consisting of braids in which the strands (or untwisted ribbons) are not permuted. The three braids $\sigma_1^2$, $\sigma_2^2$ and $\sigma_3^2$\footnote{It is wortwhile noting that these braids coincide with the generators $A^{-1}_{12}=\sigma_1^2$, $A^{-1}_{13}=\sigma_2\sigma_1^2\sigma_2^{-1}$, and $A^{-1}_{23}=\sigma_2^2$ which give a presentation of the pure braid group $PB_3$ in \cite{artin1947theory} since $\sigma_2^{-1}\sigma_1^2\sigma_2[0,0,0]=\sigma_3^2[0,0,0]$.} are in the kernel since
\begin{eqnarray}
\nonumber\sigma_1^2[a,b,c]\rightarrow [a+1,b+1,c-1],\\
\nonumber\sigma_2^2[a,b,c]\rightarrow [a-1,b+1,c+1],\\
\nonumber\sigma_3^2[a,b,c]\rightarrow [a+1,b-1,c+1].
\end{eqnarray}
It follows that, as braided 3-belts, $\sigma_1^2$, $\sigma_2^2$ and $\sigma_3^2$ commute among themselves, that is $\sigma_i^2\sigma_j^2[a,b,c]=\sigma_j^2\sigma_i^2[a,b,c]$. Any pure 3-belt $P$ can therefore be written in a standard form
\begin{eqnarray}\label{purebraid}
P=(\sigma_1^2)^{\alpha}(\sigma_2^2)^{\beta}(\sigma_3^2)^{\gamma}[0,0,0],
\end{eqnarray}
where $\alpha,\beta,\gamma$ are integers  to be determined. The length of a braid word corresponding to a pure braid is always even. Writing $P$ in pure twist form (using methods described in \cite{Bilson-Thompson2009}) $[a,b,c]$ and solving for $\alpha,\beta,\gamma$ yields (uniquely)
\begin{eqnarray}\label{eq12}
\alpha=\frac{1}{2}(a+b),\qquad \beta=\frac{1}{2}(b+c),\qquad \gamma=\frac{1}{2}(c+a),
\end{eqnarray}
which is satisfied if and only if $a,b,c$ are either all even or all odd. The word problem for pure 3-belts is therefore easily solved.  

\subsection{The braid word of a braided 3-belt}  

The image of a braid $B\in B_n$ under the map $\pi:B_n\rightarrow S_n$ is a permutation of $n$ objects. For the braid group $B_3$, the map $\pi$ therefore divides $B_3$ into six distinct cosets. An arbitrary braid in $B_n$ can then be written as $B=CP$\footnote{Because the conjugacy problem is solvable for braid groups one could write instead $B=P'C$ instead where $C$ is as before and $P'\in PB_n$ but not necessarily equal to $P$.} where $C$ is a fixed chosen coset representative and $P\in PB_n$, the kernel of $\pi$. One can then write $B$ in a canonical form by combing $P$.

Due to the lack of isogeny and subsequent group structure for braided 3-belts, extra care must be taken in determining the braid word for a general $[a,b,c]$. For a braided 3-belt with braiding $B\in B_3^c$, one can write $CB=P$ where $C$ is some permutation inducing braid in $B_3^c$ and $P\in PB_3^c$. The strategy is therefore to find a $C$ such that $C[a,b,c]$ is a pure braid ($a,b,c$ all even or all odd integers). Once $C$ is known, one can then write $B=C^{-1}P$, and determine the braid word for $P$ using equation (\ref{purebraid}).

As an example, consider $[a,b,c]=[2,4,3]$:
\begin{eqnarray}
\nonumber[2,4,3]&=&(\sigma_2\sigma_1)(\sigma_1^{-1}\sigma_2^{-1})[2,4,3],\\
\nonumber &=& (\sigma_2\sigma_1)[2,2,4],\qquad (\textrm{where}\;[2,2,4]\in PB_3^c),\\
\nonumber &=& (\sigma_2\sigma_1)(\sigma_1^2)^2(\sigma_2^2)^3(\sigma_3^2)^3,\\
\nonumber &=& \sigma_2\sigma_1^5\sigma_2^6\sigma_3^6.
\end{eqnarray}
The choice of $C$ is not unique. One could instead choose $C=\sigma_1^{-1}\sigma_3^{-1}$\footnote{This is because $\sigma_1^{-1}\sigma_2^{-1}[a,b,c]=[c-1,a,b]$ and $\sigma_1^{-1}\sigma_3^{-1}[a,b,c]=[b,c-1,a]$, and so either both are pure braids or neither are pure braids.}
\begin{eqnarray}
\nonumber[2,4,3]&=&(\sigma_3\sigma_1)^{-1}(\sigma_1^{-1}\sigma_3^{-1})[2,4,3],\\
\nonumber &=& (\sigma_3\sigma_1)[4,2,2],\\
\nonumber &=& (\sigma_3\sigma_1)(\sigma_1^2)^3(\sigma_2^2)^2(\sigma_3^2)^3,\\
\nonumber &=& \sigma_3\sigma_1^7\sigma_2^4\sigma_3^6.
\end{eqnarray}
A standard form however can defined (for example) by choosing $C$ such that is has minimal length and does not contain $\sigma_3$. 

For the case where $a,b,c$ are all integers (corresponding to an even permutation), $C$ can be chosen uniquely from the set
\begin{eqnarray}
C\in\left\lbrace 1, \sigma_1^{-1}\sigma_2^{-1}, \sigma_2^{-1}\sigma_3^{-1}, \sigma_3^{-1}\sigma_1^{-1}\right\rbrace,
\end{eqnarray}
Similarly, for the case $a,b,c$ are all half integers (corresponding to an odd permutation), $C$ can be uniquely chosen from the set 
\begin{eqnarray}
C\in \left\lbrace  \sigma_1^{-1}, \sigma_2^{-1}, \sigma_3^{-1}, \sigma_1^{-1}\sigma_2^{-1}\sigma_1^{-1}\right\rbrace. 
\end{eqnarray}
For example
\begin{eqnarray}
\nonumber\left[ \frac{7}{2},\frac{1}{2},\frac{3}{2}\right] &=&\sigma_3\sigma_{3}^{-1}\left[ \frac{7}{2},\frac{1}{2},\frac{3}{2}\right],\\
\nonumber &=& \sigma_3\left[ 1,1,3\right],\\
\nonumber &=& \sigma_3(\sigma_1^2)^1(\sigma_2^2)^2(\sigma_3^2)^2,\\
\nonumber &=& \sigma_3\sigma_1^2\sigma_2^4\sigma_3^4.
\end{eqnarray}
Similarly,
\begin{eqnarray}
\nonumber\left[ \frac{3}{2},\frac{3}{2},\frac{3}{2}\right] &=& (\sigma_1\sigma_2\sigma_1)(\sigma_1\sigma_2\sigma_1)^{-1}\left[ \frac{3}{2},\frac{3}{2},\frac{3}{2}\right],\\
\nonumber &=&(\sigma_1\sigma_2\sigma_1)\left[ 1,1,1\right],\\
\nonumber &=&(\sigma_1\sigma_2\sigma_1)\sigma_1^2\sigma_2^2\sigma_3^2.
\end{eqnarray}
Instead of $\sigma_1^{-1}\sigma_2^{-1}\sigma_1^{-1}$ one could instead choose $\sigma_2^{-1}\sigma_3^{-1}\sigma_2^{-1}$ or $\sigma_3^{-1}\sigma_1^{-1}\sigma_3^{-1}$. Using this procedure any orientable braided 3-belt may be written in braid only form in a canonical form.

\section{Jones polynomial for knotted boundaries of braided 3-belts}

Braided 3-belts are surfaces with boundary. In general, this boundary may have multiple components forming a link. When the boundary consists of a single component, this boundary will be a knot. It is not difficult to see that the boundary of a braided belt will form a single component when the pure twist vector $[a,b,c]$ consists of three half integers. This is so because in that case the edges of each ribbon are interchanged. The boundary will also form a knot in the case where the twist vector contains precisely one integer. However in this case the braided belt will not be orientable, and hence this case will not be considered. The boundary link of belts whose pure twist vectors contain more than one integer necessarily have two or three components. Three components is the maximum possible number which occurs when either $a,b,c$ are all integers. 

Alternatively, the action of a single braid generator on $[a,b,c]$ increases or decreases each $a,b,c$ by one half such that the sum $a+b+c$  increases by one half (decreases when an inverse braid generator is used). It follows that the action of an odd number of braid generators on the identity $[0,0,0]$ gives an $[a,b,c]$ with $a,b,c$ all odd multiples of $1/2$. Hence when the braid word has odd length the boundary will be a knot.

The motivation for investigating knotted boundaries is twofold. First, a knot (as well as a link) has associated with it known invariants. It may then be possible to relate these to conserved quantum numbers. Secondly, by considering which knots and links can arise as the boundaries of braided 3-belts, it may be possible to find alternative identifications of braids with SM fermions. One example is a model of leptons and quarks in terms of quantum trefoil knots, described in terms of representations of the quantum group $SU_q(2)$ by Finkelstein \cite{Finkelstein2012,Finkelstein2010,Finkelstein2007}. This model will be discussed in more detail in the next section.

\subsection{The Jones polynomial for braided belts with knotted boundary}
  
\begin{figure}[h!]
\centering
   \includegraphics[scale=0.25]{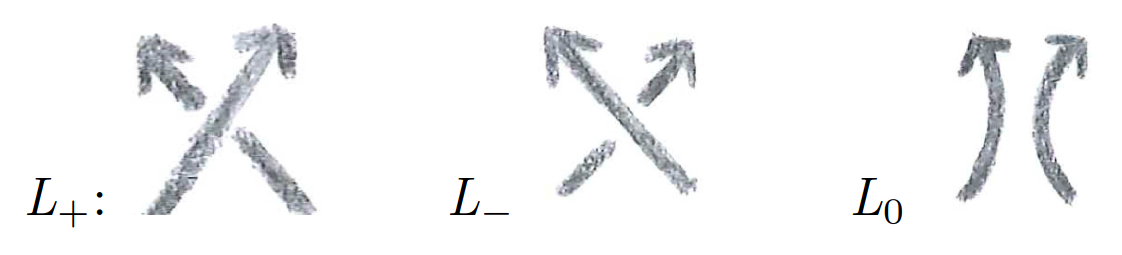}
\caption{Oriented links differing only locally.}
\label{3crb2}
\end{figure}

The skein relation for the Jones polynomial when $V(L_+)$, $V(L_-)$, and $V(L_0)$ are three oriented links differing only locally according to the diagrams in Figure \ref{3crb2},
is given by
\begin{equation}
t^{-1}V(L_+)-tV(L_-)=(t^{\frac{1}{2}}-t^{-\frac{1}{2}})V(L_0),
\end{equation}
 which can be rewritten as
\begin{equation}
V(L_+)=\alpha V(L_-)+\beta V(L_0),
\end{equation}
where $\alpha\equiv t^2$ and $\beta\equiv t^{\frac{3}{2}}-t^{\frac{1}{2}}$.

\begin{figure}[h!]
\centering
   \includegraphics[scale=0.5]{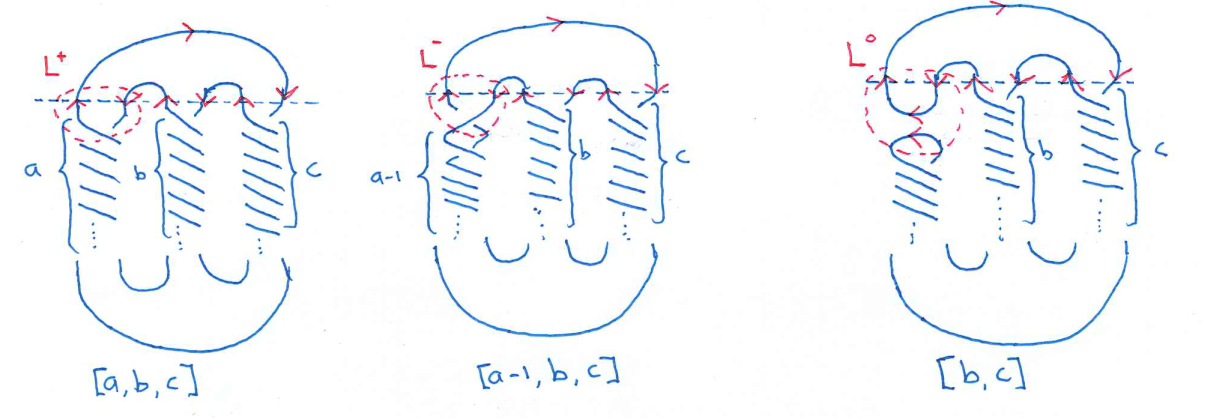}
\caption{The application of the skein relation to a braided 3-belt}
\label{3crb}
\end{figure}

Consider a braided 3-belt written in pure twist form with twist word $V(L_+)=V_{a,b,c}=[a,b,c]$ where $a,b,c$ are all positive half integers. In order to apply the skein relation one must first assign an orientation to the boundary of the belt. We will assign a clockwise direction to the boundary starting at the top of the node. It immediately follows that 
$V(L_-)=[a-1,b,c]$ and $V(L_0)=[b,c]$. That is, application of the skein relation leads to two new braided belts, one identical to the original but with one $2\pi$ twist removed from the first ribbon, and a second where the first ribbon has been removed entirely, leaving a belt composed of two ribbons. This situation is shown in Figure \ref{3crb}

By repeated application of the skein relation, the total twist on the first ribbon may be reduced to $1/2$. Next, repeated application of the skein relation to the second, and subsequently third ribbon systematically removes all but half a twist from each of the ribbons. The result of this will be a braided belt with twist word $[1/2,1/2,1/2]$ as well as various belts composed of two ribbons. For the latter the same strategy may be used to systematically remove all twists until until one is left with the twist word $[1/2,1/2]$. Now $[1/2,1/2,1/2]$ and $[1/2,1/2]$ are identified with the trefoil knot and Hopf link respectively, whose Jones polynomials are known. 

With this strategy in mind, and with $V_{a,b,c}$ representing the Jones polynomial of the boundary of the braided belt with pure twist word $[a,b,c]$, we write,
\begin{equation}
V_{a,b,c}=\alpha V_{a-1,b,c}+\beta V_{b,c},
\end{equation}
which after iterating $a-\frac{1}{2}$ times gives
\begin{equation}
V_{a,b,c}=\alpha^{a-\frac{1}{2}}V_{\frac{1}{2},b,c}+\beta\frac{\alpha^{a-\frac{1}{2}}-1}{\alpha-1}V_{b,c}.
\end{equation}
Repeatedly applying the skein relation $b-\frac{1}{2}$ and $c-\frac{1}{2}$ times to the second and third ribbon respectively and rearranging gives
\begin{eqnarray}
\nonumber J_{a,b,c}=\alpha^{a+b+c-\frac{3}{2}}V_{\frac{1}{2},\frac{1}{2},\frac{1}{2}}+\beta\frac{\alpha^{a+b+c-\frac{3}{2}}-\alpha^{a+b-1}}{\alpha-1}V_{\frac{1}{2},\frac{1}{2}}+\beta\frac{\alpha^{a+b-1}-\alpha^{a-\frac{1}{2}}}{\alpha-1}V_{\frac{1}{2},c}\\+\beta\frac{\alpha^{a-\frac{1}{2}}-1}{\alpha-1}V_{b,c}.
\end{eqnarray}

Applying the same strategy to the two-ribbon braids $V_{\frac{1}{2}, c}$ and $V_{b,c}$ and using the fact that $V_{c}=V(\textrm{unknot})=1$, we obtain
\begin{eqnarray}
V_{b,c}=\alpha^{b-\frac{1}{2}}V_{\frac{1}{2},c}+\beta\frac{\alpha^{b-\frac{1}{2}}-1}{\alpha-1},\\
V_{\frac{1}{2},c}=\alpha^{c-\frac{1}{2}}V_{\frac{1}{2},\frac{1}{2}}+\beta\frac{\alpha^{c-\frac{1}{2}}-1}{\alpha-1},
\end{eqnarray}
and hence
\begin{eqnarray}
\nonumber V_{a,b,c}&=&\alpha^{a+b+c-\frac{3}{2}}V_{\frac{1}{2},\frac{1}{2},\frac{1}{2}}+\frac{\beta}{\alpha-1}(3\alpha^{a+b+c-\frac{3}{2}}-\alpha^{a+b-1}-\alpha^{a+c-1}-\alpha^{b+c-1})V_{\frac{1}{2},\frac{1}{2}} \\
&+&\frac{\beta^2}{(\alpha-1)^2}(3\alpha^{a+b+c-\frac{3}{2}}-\alpha^{a+c-1}-\alpha^{b+c-1}-\alpha^{a+b-1}).
\end{eqnarray}
The above formula expresses the Jones polynomial for an arbitrary twist word $[a,b,c]$ with half integers $a,b,c$ in terms of the Jones polynomial of the trefoil knot ($V(\mathrm{Trefoil})=V_{\frac{1}{2},\frac{1}{2},\frac{1}{2}}=t+t^3-t^4$) and the Hopf link ($V(\mathrm{Hopf Link})=V_{\frac{1}{2},\frac{1}{2}}=-(t^{\frac{5}{2}}+t^{\frac{1}{2}})$). The above expression is fully symmetric in $a,b,c$ as expected. Note that although in the above derivation of the Jones polynomial it was assumed explicitly that $a,b,c$ are positive half integers, the formula works also for negative half integers or a mixture of positive and negative half integers.

\subsection{Knotted boundaries of minimally braided belts}

We now use the Jones polynomial formula to investigate the simplest non-trivially knotted boundaries. By simplest we mean the focus is on the shortest braid words whose boundaries correspond to non-trivial knots. The simplest case occurs when the length of the braid word is equal to three. The first step is to convert the braid word into a pure twist vector so that the Jones polynomial formula may be used. 

Because the formula for the Jones polynomial is symmetric in $a,b,c$ (reflecting the fact that the order of $a,b,c$ in the twist vector is irrelevant), together with the braid relations of the braid group, the topologically distinct boundaries that are possible is greatly reduced. Converting a braid word into a pure twist vector, every braid generators increases or decreases each $a,b,c$ by one half in such a way that the total sum $a+b+c$ increases by one half for each braid generator and decreases by one half for each inverse generator. For a braid word of length three this means that $a,b,c\in \left\lbrace -\frac{3}{2},-\frac{1}{2},\frac{1}{2},\frac{3}{2} \right\rbrace$, and also $a+b+c\in \left\lbrace -\frac{3}{2},-\frac{1}{2},\frac{1}{2},\frac{3}{2} \right\rbrace$.

Consider first the case where $a+b+c=\frac{3}{2}$, which consists of combinations of three braid generators and no inverse braid generators. The possibilities for $a$, $b$, and $c$ (remembering the symmetries between them) are given in Table \ref{3generatorsa}.
\begin{table}[h!]
  \begin{center}
\begin{tabular}{|c|c|c|c|}
\hline 
a & b & c & knot \\ 
\hline 
3/2 & 3/2 & -3/2 & $9_{46}$ \\ 
\hline 
3/2 & 1/2 & -1/2 & unknot\\ 
\hline
1/2 & 1/2 & 1/2 & $3_1$ \\ 
\hline 
\end{tabular}
\caption{Knotted boundaries of a braided belt composed of three braid generators.}
\label{3generatorsa}
\end{center}
\end{table}
Interestingly, each of the three possibilities gives a different knot. The labeling of knots follows the Alexander-Briggs notation as found in \cite{knotatlas}.

Next consider the case where $a+b+c=\frac{1}{2}$ (these braid words contain two braid generators and one inverse braid generator). The possibilities are listed in Table \ref{3generatorsb}.  
\begin{table}[h!]
  \begin{center}
\begin{tabular}{|c|c|c|c|}
\hline 
a & b & c & knot\\ 
\hline 
3/2 & 1/2 & -3/2 & $6_1$ \\ 
\hline 
3/2 & -1/2 & -1/2 & $4_1$ \\ 
\hline
1/2 & 1/2 & -1/2 & unknot\\ 
\hline 
\end{tabular}
\caption{Knotted boundaries of a braided belt composed of two braid generators and one inverse braid generator.}
\label{3generatorsb}
\end{center}
\end{table}
The remaining two cases ($a+b+c=-\frac{1}{2}$ and $a+b+c=-\frac{3}{2}$) give simply the mirror images of the knots already obtained. With the exception of the knot $4_1$, the mirror images are distinct knots which we label using as asterisk, for example $3_1^*$. 

More generally one may investigate longer braids. For the case where the braid word length is an odd integer $N$ (required for the boundary to form a knot)
\begin{eqnarray}
\nonumber a,b,c\in \left\lbrace -\frac{N}{2}, -\frac{N-2}{2},...,\frac{N-2}{2},\frac{N}{2}\right\rbrace,\\
\nonumber a+b+c\in \left\lbrace -\frac{N}{2}, -\frac{N-2}{2},...,\frac{N-2}{2},\frac{N}{2}\right\rbrace.
\end{eqnarray}

In the above analysis, only in one case do $a,b,c$ all have the same sign. This is when $a=b=c=\frac{1}{2}$ (respectively $a=b=c=-\frac{1}{2}$ for the inverse) with the boundary corresponding to a trefoil knot. The braid word in this case is $\sigma_1\sigma_2\sigma_1$ (and its inverse). One of the assumptions of the Helon model is the no charge mixing assumption which states that a braided belt is not allowed to mix clockwise and counterclockwise twisting on its ribbons. This assumption is simply a restatement in a different context of a similar assumption made in the models of Harari and Shupe \cite{Harari1979,Shupe1979}. Imposing a similar restriction on pure twist vectors for the simplest case of word length three then singles out the trefoil knots as the only possible knotted boundaries. Not only are the trefoil knots the simplest knots, they also form the basis of an alternative model of leptons and quarks proposed by Finkelstein. It may then be possible to use the present context of braided belts to embed that model, making it compatible with LQG. This possibility is discussed in the next section.

Finally, in each case above, the boundary knot always corresponds to a prime knot. This remains the case at least up to $a+b+c=9/2$. However a more complete study needs to be carried out. This is currently work in progress.

\section{Relation to quantum trefoil knots and dynamics}

In that model, electric charge is related to the writhe and rotation of the knots. Gauge bosons are formed via the knot sum of two trefoil knots. Through a better understanding of the boundaries of braided belts it may be possible to in the future to find an embedding of Finkelstein's quantum trefoil knots into ribbon networks.  

\section{A trefoil knot model of leptons and quarks based on representations of $SU_q(2)$}

The Helon model represents leptons and quarks as certain braided belts which may be embedded within trivalent and tetravalent spin networks. An alternative model proposed by Finkelstein instead represents these fundamental particles as oriented trefoil knots labeled by irreducible representations of the quantum group $SU_q(2)$ \cite{Finkelstein2012,Finkelstein2010,Finkelstein2007}. Each of the four classes of fermions is associated with one of the four oriented trefoil knots. 

For an oriented knot, each vertex/crossing can be identifies as either being positive one or negative one, depending on whether the orientation of the overline is carried into the orientation of the underline by a counterclockwise or clockwise rotation respectively. The sum of all the crossings signs is the writhe, $w$, whereas the rotation $r$ of a knot is simply the number of rotations of the tangent in going one around the knot. Together with the crossing number $N=3$, the ooriented trefoil knots representing families of fermions are then identified with the $j=3/2$ representations of the quantum group $SU_q(2)$, shown in Table \ref{qreps}.

\begin{table}[h!]
\renewcommand{\arraystretch}{2.0}
  \begin{center}
\begin{tabular}{|c|c|c|c|}
\hline 
$(w,r)$ & $D^{3/2}_{\frac{w}{2}\frac{r+1}{2}}$ & $Q$ & Family \\
\hline
(-3,2) & $D^{3/2}_{-\frac{3}{2}\frac{3}{2}}$ & 0 & $(\nu_e,\nu_{\mu},\nu_{\tau})$ \\ 
\hline 
(3,2) & $D^{3/2}_{\frac{3}{2}\frac{3}{2}}$ & -1 & $(e,\mu,\tau)$\\ 
\hline
(3,-2) & $D^{3/2}_{\frac{3}{2} -\frac{1}{2}}$ & $-\frac{1}{3}$ & $(d,s,b)$ \\ 
\hline
(-3,-2) & $D^{3/2}_{-\frac{3}{2} \frac{1}{2}}$ & $\frac{2}{3}$ & $(u,s,t)$ \\ 
\hline 
\end{tabular}
\caption{The four families of fermions as $j=3/2$ representations of the quantum group $SU_q(2)$.}
\label{qreps}
\end{center}
\end{table}

\begin{figure}[h!]
\centering
   \includegraphics[scale=0.7]{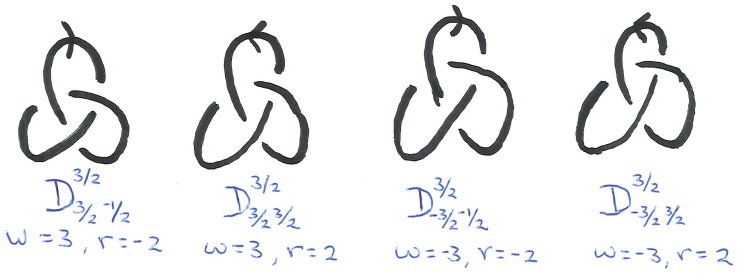}
\caption{Trefoil knots as representations of $SU_q(2)$.}
\label{trefoilfig}
\end{figure}

The representations of $SU_q(2)$, $D^j_{mm'}$ can now be interpreted as the kinematical description of a quantum state according to the relation
\begin{eqnarray}
(j,m,m')=\frac{1}{2}(N,w,\pm r+1),
\end{eqnarray}
where $(N,w,r)$ describes a classical knot. This equation establishes a connection between a quantized knot described by $D^{N/2}_{\frac{w}{2}\frac{\pm r+1}{2}}$ and a classical knot described by $(N,w,r)$. For the trefoil knots this correspondence is one-to one as illustrated in Figure \ref{trefoilfig}, although in general it is not. The electric charge of a fermion is determined by the formula
\begin{eqnarray}
Q=-\frac{1}{3}(m+m').
\end{eqnarray}

It the previous section it was shown that the trefoil knots can arise as the boundaries of braided belts. This might makes it possible to cast Finkelstein's model into the framework of braided 3-belts discussed here. As in the Helon model, leptons and quarks in Finkelstein's knot model might then be considered as topological defects in a background quantum spacetime provided by LQG. Such an embedding of the knot model into quantum gravity has not been considered before. 

One of the outstanding open probelms in the Helon model is developing the dynamics. In \cite{wan2009effective} it is suggested that in a dynamical theory of interacting braids, gauge symmetries may be incorporated by including the spin network labels. An alternative approach may be to instead look at braided belts with knotted boundaries and associating these boundaries with irreducible representations of $SU_q(2)$ as in the knot model of Finkelstein. For example, to develop the electroweak dynamics in the knot model, one retains the formal structure of standard electroweak but at the same time expands the quantum fields involved in a new class of normal modes that are determined by the algebra of $SU_q(2)$ \cite{finkelstein2007field}.  The knot Lagrangian is constructed by replacing the left chiral field operators for the elementary fermions, $\Psi_L(x)$, and the electroweak vectors $W_{\mu}(x)$, of the SM by the knotted field operators $\Psi^{\frac{3}{2}}_L(x)D^{\frac{3}{2}}_{mm'}$ and $W^3_{\mu}(x)D^3_{mm'}$ respectively. This Lagrangian preserves the local $SU(3)\times SU(2)\times U(1)$ and thus the dynamics. The knot factors however introduce form factors for all terms. Connecting the knot model with the Helon model by identifying the knots as the boundaries of braided 3-belts then provides a unique approach to developing a dynamical model of interacting braids. The details remain to be worked out however. 

\section{Discussion}

This paper has carried on earlier work on braided belts in which several approaches to manipulating braided belts were considered \cite{Bilson-Thompson2009}. Whereas in that paper it was shown that a braided 3-belt may always be written in pure twist form with trivial braiding, in the first part of this paper it was shown that any orientable 3-belt may also be written in a braid only form with trivial twisting. An algorithm for determining the braid word was presented, allowing any (orientable) braided 3-belt to be written in a canonical form. 

When written in braid only form a connection between the Helon model and a complementary descriptions of leptons and quarks in terms of the ideals of normed division algebras becomes apparent \cite{furey2016standard,gresnigt2017braidsgroups,gresnigt2018braids}. Specifically, writing the Helon braids in braid only form, and suitably identifying products of pairs of braid generators with raising and lowering operators, the Helon braids coincide precisely with the basis states of the complex octonions (or equivalently, the Clifford algebra $C\ell(6)$ \cite{stoica2017standard}).

Associated with any braided belt is its boundary link, and the reduced boundary link. The latter is invariant under the standard Pachner moves which generate the evolution algebra of braided belts embedded in trivalent ribbon networks. In the second half of this paper it was determined when the boundary link of an orientable braided belt corresponds to a knot. This is the case when the pure twist vector $[a,b,c]$ consists of three half integers. For this case a formula for the Jones Polynomial was determined.

Considering knotted boundaries is useful for a number of reasons. Because the trefoil knots arise as boundaries of braided belts (in this case with pure twist vectors $[\frac{1}{2},\frac{1}{2},\frac{1}{2}]$) it is possible to embed an independent model of SM fermions based on quantum trefoil knots into trivalent and tetravalent spin networks. At the same time, it shows that the Helon model mapping from braided belts to leptons and quarks need not necessarily be unique. By considering the writhe, rotation, and crossing number of the boundary knots instead of braiding and twisting, instead of the braiding and twisting of the 3-belt itself,  additional representations of SM fermions are possible. The knot model of Finkelstein provides the simplest example of this. There likely exist some relationships between the braiding and twisting of braided 3-belts with the writhe and rotation of its boundary although what these relationships are remains to be determined.

It is also possible to associate particular representations of the quantum group $SU_q(2)$ to the boundaries of braided belts. The quantum trefoil knots in Finkelsteins model correspond to the $j=3/2$ representations. This might provide a novel means of developing the gauge symmetries of interacting braids in the future, thereby extending earlier work on the dynamics by Smolin and Wan \cite{smolin2008propagation} and Wan \cite{wan2009effective}. In the latter it was suggested that gauge symmetry can be developed by including spin network labels. The work in the present paper suggests a complementary approach. The details of both approaches remain to be worked out. 

A deeper connection between the knot model of Finkelstein and the Helon model may be possible. In the former, fermions correspond to the $j=3/2$ representations of $SU_q(2)$. One can extend the model into a preon model consistent with the Harari Shupe model by identifying the fundamental $j=1/2$ representations with preons. Since both the knot and Helon model reproduce the Harari Shupe model, it may be possible to identify the fundamental representations of $SU_q(2)$ with helons (a twisted ribbon). This requires further investigation. 

It would also be interesting to investigate more carefully the spectrum of knots that may arise as the boundary of braided belts. For the cases considered in this paper a knotted boundary always corresponds to a prime knot. It is not known if this is true in general however. This is currently under investigation.

\section{Acknowledgments}

This work is supported by the National Natural Science Foundation of China grant RRSC0116. The author wishes to thank Adam Gillard, and Benjamin Martin for insightful discussions and Benjamin Martin for hosting the author at the University of Aberdeen in July 2017.
\bibliography{NielsReferences}  
\bibliographystyle{unsrt}  

\end{document}